\documentclass[12pt]{amsart}
\usepackage[english]{babel}

\title{A combinatorial theorem of Schur implies infinitude of  primes} 

\author{Labib Haddad}
\address{120 rue de Charonne, 75011 Paris, France}
\email{labib.haddad@wanadoo.fr}

\usepackage{amssymb}
\usepackage{amsmath}

\newcommand{\su}{\subsection*}
\newcommand{\head}{\section*}
\newcommand{\noi}{\noindent}

\newcommand{\leqs}{\leqslant}

\newcommand{\geqs}{\geqslant}

\newcommand{\lopar}{\noi \{$\looparrowright$ \ }

\begin{document}
\maketitle

\thispagestyle{empty}

\markboth{Labib Haddad}{Infinitude of primes}

\

\hfill {\it Much ado about ...}

\hfill {William Shakespeare}

\

\

In a paper, [1], recently published,  Elsholz shows \lq\lq that Fermat's last theorem and a combinatorial theorem of Schur on monochromatic solutions of $a + b = c$ implies that there exist infinitely many primes."

The proof can be greatly simplified, making no use at all of Fermat's last theorem,  and using only a weak form of the theorem of Schur.

Of course the reference to Fermat's last theorem sounds more attractive, somehow, but it is really not necessary, it can be dispensed with.

\

Let us now further elaborate !

\

\head{The theorem of Schur}

\

This result is usually stated as follows.

\su{Theorem of Schur} Let $t>0$ be any integer. For some large enough integer $s$, if the integers in the interval $[1,s]$ are  colored using $t$ distinct colors, there is (at least) one monochromatic  solution to $a+b =c$ such that $1 \leqs a <  b < c \leqs s$.

\

Clearly, the following  weaker version of the theorem also holds.

\su{Schur's lemma} If  the positive integers are distributed in a finite number of pigeon-holes, then one of the pigeon-holes (at least) contains three integers, $a < b < c$, such that $a+b = c$.

\

Notice that Schur's lemma acts like a {\bf selector}: It says that, in any situation where you classify the positive integers into a finite number of classes,  a special solution to the equation $a+b= c$ can be {\bf selected}, the three integers in the same class !

\

\head{If the primes were finite in number}

\

Suppose there are only $k$ primes, say, $p_1 < p_2 < \dots < p_k$. Much the same as in the proof of Elsholz,  notice that each  positive integer $m$ can be written as
$$m = p_1^{4q_1 + r_1} \dots p_k^{4q_k + r_k} \ , \ q_i \geqs 0 \ ,\ 0\leqs r_i < 4.$$
Setting
$$u(m) = p_1^{q_1} \dots p_k^{q_k} \ , \ M(m) = p_1^{r_1}\dots p_k^{r_k},$$
obtains
$$m = u(m)^4 M(m).$$
As Elsholz observes, there is no need to know that this representation of the integer $m$ is unique. 

\

\noi Anyway, call the integer $M(m)$ {\bf the mantissa} of $m$, for brevity's sake. There are only $4^k$ mantissas. Distribute the positive integers in $4^k$ pigeon-holes, each containing integers $m$ with the same mantissa!

\

\noi Let $a < b < c$ be a selected solution  to the equation $a+b = c$, by Schur's lemma.  So,
$a,b,c,$ have a same mantissa, $d$.  Thus, we have
$$a = u(a)^4d \ , \ b=u(b)^4 d\ , \ c = u(c)^4d,$$
$$u(a)^4 d+ u(b)^4d= u(c)^4d,$$
$$u(a)^4 + u(b)^4= u(c)^4$$
$$u(c)^4 - u(b)^4 = (u(a)^2)^2.$$
But there is no positive integers $x,y,z,$ such that $z^4 - y^4 = x^2$, by {\it infinite descent}. A contradiction! The primes are infinite in number!\qed

\

\lopar For that matter, the whole thing can be dismissed arguing as follows: Let $E$ be the statement that there is only a finite number of primes. If statement $E$ were true, then Arithmetics would be an inconsistent theory, so that its theorems, all, would be true and false ! So, any of the theorems of Arithmetics implies the infinitude of primes!

Yes, this looks a bit {\it specious}! Yet, the argument is {\it tight} enough! Indisputable!\}

\

\

\

\

\

\head{REFERENCES}

\

\noi [1] Christian Elsholtz,  {\sl Fermat's last theorem implies Euclid's infinitude of primes}, Amer. Math. Monthly, {\bf 128} (2021), no 3, 250-257.

\

\

\enddocument